\documentclass[12pt,draft]{amsart}
\usepackage[all]{xy}
\usepackage{amsfonts}
\usepackage{amssymb}
\usepackage{enumerate}

\usepackage{color}

\newtheorem{teo}{Theorem}[section]
\newtheorem{lem}[teo]{Lemma}

\newtheorem{cor}[teo]{Corollary}

\theoremstyle{definition}
\newtheorem{dfn}[teo]{Definition}
\newtheorem{rk}[teo]{Remark}
\newtheorem{ex}[teo]{Example}

\def\<{\langle}
\def\>{\rangle}

\def\b{\beta}

\def\l{{\lambda}}
\def\r{\rho}

\def\t{\tau}

\def\G{{\Gamma}}

\def\C{{\mathbb C}}

\def\Z{{\mathbb Z}}

\def\End{\mathop{\rm End}\nolimits}
\def\Ker{\mathop{\rm Ker}\nolimits}
\def\Im{\mathop{\rm Im}\nolimits}
\def\Aut{\operatorname{Aut}}

\def\Id{\operatorname{Id}}

\def\1{\mathbf 1}

\newcommand{\ov}[1]{\overline{#1}}
\newcommand{\til}[1]{\widetilde{#1}}
\newcommand{\wh}[1]{\widehat{#1}}

\begin{document}

\title[Twisted conjugacy separable groups]
{Twisted conjugacy separable groups}

\author{Alexander Fel'shtyn}
\address{Instytut Matematyki, Uniwersytet Szczecinski,
ul. Wielkopolska 15, 70-451 Szczecin, Poland and Department of Mathematics,
Boise State University,
1910 University Drive, Boise, Idaho, 83725-155, USA}
\email{felshtyn@diamond.boisestate.edu, felshtyn@mpim-bonn.mpg.de}

\author{Evgenij Troitsky}
\thanks{The second author is partially supported by
RFFI Grant  05-01-00923
and Grant ``Universities of Russia''}
\address{Dept. of Mech. and Math., Moscow State University,
119992 GSP-2  Moscow, Russia}
\email{troitsky@mech.math.msu.su}
\urladdr{
http://mech.math.msu.su/\~{}troitsky}


\begin{abstract}
We study the notion of twisted conjugacy separability
(essentially introduced in our previous
paper for a proof of twisted version of Burnside-Frobenius theorem)
and some related properties.
We give examples of groups with and without this property and study
its behavior under some extensions.
An affirmative
answer to the twisted Dehn conjugacy problem
for polycyclic-by-finite group is obtained.
Some problems for the further study
are indicated.
\end{abstract}

\maketitle

\tableofcontents

\section{Introduction}

\begin{dfn}
Let $G$ be a countable discrete group and $\phi: G\rightarrow G$ an
endomorphism.
Two elements $x,x'\in G$ are said to be
 $\phi$-{\em conjugate} or {\em twisted conjugate,}
iff there exists $g \in G$ with
$
x'=g  x   \phi(g^{-1}).
$
We shall write $\{x\}_\phi$ for the $\phi$-{\em conjugacy} or
{\em twisted conjugacy} class
 of the element $x\in G$.
The number of $\phi$-conjugacy classes is called the {\em Reidemeister number}
of an  endomorphism $\phi$ and is  denoted by $R(\phi)$.
If $\phi$ is the identity map then the $\phi$-conjugacy classes are the usual
conjugacy classes in the group $G$.
\end{dfn}

If $G$ is a finite group, then the classical Burnside{-Frobenius}
theorem (see e.g. \cite{serrerepr},
\cite[p.~140]{Kirillov})
says that the number of
classes of irreducible representations is equal to the number of conjugacy
classes of elements of $G$.  Let $\wh G$ be the {\em unitary dual} of $G$,
i.e. the set of equivalence classes of unitary irreducible
representations of $G$.
The attempts to generalize this theorem to the case of non-identical
automorphism and of non-finite group
(i.e., to identify the Reidemeister number of $\phi$ and the number
of fixed points of $\wh\phi$ on an appropriate dual object of $G$,
provided that one of these numbers is finite)
were inspired by the dynamical
questions and were the subject  of a series of papers
\cite{FelHill,FelBanach,FelshB,FelTro,FelTroVer,FelIndTro,polyc}.

In the present paper we study the following property for a countable
discrete group $G$ and its automorphism $\phi$: we say that
the group is $\phi$-\emph{conjugacy separable} if its Reidemeister
classes can be distinguished by homomorphisms onto finite groups, and
we say that it is \emph{twisted conjugacy separable} if it is
$\phi$-conjugacy separable for any automorphism $\phi$
with $R(\phi)<\infty$ (\emph{strongly twisted conjugacy separable}, if
we remove this finiteness restriction)
(Definitions \ref{dfn:ficonjsep} and \ref{dfn:twisconjsep}).
This notion was used in \cite{polyc} to prove the twisted Burnside-Frobenius
theorem for polycyclic-by-finite groups with the finite-dimensional part
of the unitary dual $\wh G$ as an appropriate dual object.
Related questions were studied in \cite{BardBokutVesnin}, \cite{Martino}.

After some preliminary considerations we prove \textbf{the main
results} of the paper, namely

\begin{enumerate}
    \item {\bf Classes of twisted conjugacy separable groups:} Polycyclic-by-finite groups
    are strongly twisted conjugacy separable groups
    (Theorem \ref{teo:almostpolaresep}).
    \item {\bf Twisted conjugacy separability respects some extensions:}
Suppose, there is an extension $H\to G\to G/H$, where the group $H$ is a
finitely generated characteristic
twisted conjugacy separable group; $G/H$ is finitely generated {\rm FC}-group (i.e.,
a group with finite conjugacy classes).
Then $G$ is a twisted conjugacy separable group (a reformulation of
Theorem \ref{teo:phiconjandexten}).
    \item {\bf Examples of groups, which are not twisted conjugacy separable:}
HNN, Ivanov and Osin groups (Section \ref{sec:osin}).
\item {\bf The affirmative answer to the twisted Dehn conjugacy problem for
polycyclic-by-finite groups} (Section \ref{sec:desisionprob}).
\item
{\bf Residually finite finitely generated groups are twisted conjugacy separable,}
in particular {\bf twisted Burnside-Frobenius theorem is true for them}
in the following formulation: Let $G$ be a finitely generated residually finite group and
$\phi$ its automorphism with $R(\phi)<\infty$. Then $R(\phi)=S_f(\phi)$,
where $S_f(\phi)$ is the number of fixed points of $\wh\phi:\wh G_f\to \wh G_f$,
$\wh\phi (\r)=\r\circ\phi$, where $\wh G_f$ is the part of the unitary dual
$\wh G$, which is formed by the finite-dimensional representations
(Section \ref{sec:resfin}).
\end{enumerate}

A number of examples
of groups and automorphisms with finite Reidemeister numbers was
obtained and studied in \cite{FelshB,gowon,FelHillWong,FelTroVer,FelIndTro}.

\medskip
The interest in twisted conjugacy relations has its origins, in particular,
in the Nielsen-Reidemeister fixed point theory (see, e.g. \cite{Jiang,FelshB}),
in Selberg theory (see, eg. \cite{Shokra,Arthur}),
and  Algebraic Geometry (see, e.g. \cite{Groth}).
The congruences give some necessary conditions for the realization problem
for Reidemeister numbers in topological dynamics.

\medskip\noindent
{\bf Acknowledgement.} The present research is a
part of our joint research programm
in Max-Planck-Institut f\"ur Mathematik (MPI)
in Bonn.
We would like to thank the MPI for its kind support and
hospitality while the most part of this work has been completed.

The authors are grateful to
F.~Grunewald,
S.~Ivanov,
M.~Kapovich,
D.~Osin,
V.~Re\-mes\-len\-ni\-kov,
M.~Sapir,
and A.~Shtern
for helpful discussions.

\section{Preliminary considerations}\label{sec:prelim}
%%%%%%%%%%%%%%%%%%%%%%%%%%%%%

The following construction relates $\phi$-conjugacy classes and
some conjugacy classes of another group. It was  obtained
in topological context by Boju Jiang and Laixiang
Sun in \cite{Jiang88}.
Consider the action of $\Z$ on $G$, i.e. a homomorphism $\Z\to\Aut(G)$,
$n\mapsto\phi^n$. Let $\G$ be a corresponding semi-direct product
$\G=G\rtimes\Z$:
\begin{equation}\label{eq:defgamma}
    \G:=< G,t\: |\: tgt^{-1}=\phi(g) >
\end{equation}
in terms of generators and relations, where $t$ is a generator of $\Z$.
The group $G$ is a normal subgroup
of $\G$. As a set, $\G$ has the form
\begin{equation}\label{eq:razlgamm}
    \G=\sqcup_{n\in\Z}G\cdot t^n,
\end{equation}
where $G\cdot t^n$ is the coset by $G$ containing $t^n$.

\begin{rk}
Any usual conjugacy class of $\G$ is contained in some $G\cdot t^n$.
Indeed, $g g't^n g^{-1}=gg'\phi^n(g^{-1}) t^n$ and $t g't^n t^{-1}=\phi(g')t^n$.
\end{rk}

\begin{lem}\label{lem:HillRanbij}
Two elements $x,y$ of $G$ are $\phi$-conjugate iff $xt$ and $yt$ are conjugate
in the usual sense in $\G$. Therefore $g\mapsto g\cdot t$ is a bijection
from the set of $\phi$-conjugacy classes of $G$ onto the set of conjugacy
classes of $\G$ contained in $G\cdot t$.
\end{lem}

\begin{proof}
If $x$ and $y$ are $\phi$-conjugate then there is a $g\in G$ such that
$gx=y\phi(g)$. This implies $gx=ytgt^{-1}$ and therefore $g(xt)=(yt)g$
so $xt$ and $yt$ are conjugate in the usual sense in $\G$. Conversely,
suppose $xt$ and $yt$ are conjugate in $\G$. Then there is a $gt^n\in\G$
with $gt^nxt=ytgt^n$. From the relation $txt^{-1}=\phi(x)$ we obtain
$g\phi^n(x)t^{n+1}=y\phi(g)t^{n+1}$ and therefore $g\phi^n(x)=y\phi(g)$.
Hence, $y$ and $\phi^n(x)$ are $\phi$-conjugate. Thus,
$y$ and $x$ are $\phi$-conjugate, because $x$ and $\phi(x)$ are always
$\phi$-conjugate: $\phi(x)=x^{-1} x \phi(x)$.
\end{proof}

\section{Twisted conjugacy separability}\label{sec:separat}
We would like to give a generalization of the following well known notion.

\begin{dfn}
A group $G$ is \emph{conjugacy separable} if any pair $g$, $h$ of
   non-conjugate elements of $G$ are non-conjugate in some finite
   quotient of $G$.
\end{dfn}

It was proved that polycyclic-by-finite groups
are conjugacy separated (\cite{Remes,Form}, see also \cite[Ch.~4]{DSegalPoly}).
Also, residually finite recursively presented
Burnside $p$-groups constructed by R. I. Grigorchuk
\cite{GrFA} and by N. Gupta and S. Sidki \cite{GuSi}
are shown to be conjugacy separable
when $p$ is an odd prime in \cite{WiZa}.

We can introduce the following notion, which coincides with the previous
definition in the case $\phi=\Id$.

\begin{dfn}\label{dfn:ficonjsep}
A group $G$ is \emph{$\phi$-conjugacy separable} with respect to
an automorphism $\phi:G\to G$ if any pair $g$, $h$ of
   non-$\phi$-conjugate elements of $G$ are non-$\ov\phi$-conjugate in some finite
   quotient of $G$ respecting $\phi$.
\end{dfn}

This notion is closely related to the notion ${\rm RP}(\phi)$
introduced in \cite{polyc}.

\begin{dfn}\label{dfn:Rperiodi}
We say that a group $G$ has the property   {\sc RP} if for
any automorphism $\phi$ with $R(\phi)<\infty$ the
characteristic functions $f$ of {\sc Reidemeister} classes (hence
all $\phi$-central functions) are   {\sc periodic} in the
following sense.

There exists a finite group $K$,
its automorphism $\phi_K$, and epimorphism $F:G\to K$ such that
\begin{enumerate}
    \item The diagram
    $$
    \xymatrix{
G\ar[r]^\phi\ar[d]_F& G\ar[d]^F\\
K\ar[r]^{\phi_K}& K
    }
    $$
    commutes.
    \item $f=F^*f_K$, where $f_K$ is a characteristic function
    of a subset of $K$.
\end{enumerate}

If this property holds for a concrete automorphism $\phi$, we
will denote this by RP($\phi$).
\end{dfn}

One gets immediately the following statement.

\begin{teo}\label{teo:phiconjRP}
Suppose, $R(\phi)<\infty$. Then $G$ is $\phi$-conjugacy separable if and
only if $G$ is {\rm RP}$(\phi)$.
\end{teo}

\begin{proof}
Indeed, let $F_{ij}:G\to K_{ij}$ distinguish $i$th and $j$th $\phi$-conjugacy
classes, where $K_{ij}$ are finite groups, $i,j=1,\dots,R(\phi)$. Let
$F:G\to \oplus_{i,j}K_{ij}$, $F(g)=\sum_{i,j}F_{ij}(g)$, be the diagonal
mapping and $K$ its image. Then $F:G\to K$ gives RP$(\phi)$.

The opposite implication is evident.
\end{proof}

\begin{dfn}\label{dfn:twisconjsep}
A group $G$ is \emph{twisted conjugacy separable} if it is $\phi$-conjugacy
separable for any $\phi$ with $R(\phi)<\infty$.

A group $G$ is \emph{strongly twisted conjugacy separable} if it is $\phi$-conjugacy
separable for any $\phi$.
\end{dfn}

From Theorem \ref{teo:phiconjRP} one immediately obtains

\begin{cor}\label{cor:twistconjRP}
A group $G$ is twisted conjugacy separable if and only if it is {\rm RP}.
\end{cor}

\begin{teo}\label{teo:razdelgamg}
Let $F:\G\to K$ be a morphism onto a finite group $K$ which separates
two conjugacy classes of $\G$ in $G\cdot t$. Then the restriction
$F_G:=F|_G:G\to \Im(F|_G)$
separates the corresponding $($by the bijection from
Lemma {\rm \ref{lem:HillRanbij})} $\phi$-conjugacy classes in $G$.
\end{teo}

\begin{proof}
First of all let us remark that $\Ker(F_G)$ is $\phi$-invariant.
Indeed, suppose $F_G(g)=F(g)=e$. Then
$$
F_G(\phi(g))=F(\phi(g))=F(tgt^{-1})=F(t)F(t)^{-1}=e
$$
(the kernel of $F$ is a normal subgroup).

Let $gt$ and $\til gt$ be some representatives of the mentioned
conjugacy classes. Then
$$
F((ht^n)gt (ht^n)^{-1})\ne F(\til gt),\qquad\forall h\in G,\: n\in \Z,
$$
$$
F(ht^ngt )\ne F(\til gt ht^n),\qquad\forall h\in G,\: n\in \Z,
$$
$$
F(h\phi^n(g)t^{n+1} )\ne F(\til g \phi(h)t^{n+1}),\qquad\forall h\in G,\: n\in \Z,
$$
$$
F(h\phi^n(g))\ne F(\til g \phi(h)),\qquad\forall h\in G,\: n\in \Z,
$$
in particular, $F(hg\phi(h^{-1}))\ne F(\til g )$ $\forall h\in G$.
\end{proof}

\begin{teo}\label{teo:conjsepandRP}
Let some class of conjugacy separable groups be closed under
taking semidirect products by $\Z$. Then this class consists of
strongly twisted conjugacy separable groups.
\end{teo}

\begin{proof}
This follows immediately from Theorem \ref{teo:razdelgamg} and
Theorem \ref{teo:phiconjRP}.
\end{proof}

\section{First examples: polycyclic-by-finite groups}\label{sec:almostpolsep}

As an application we obtain another proof of the main theorem for
polycyclic-by-finite groups.

Let $G'=[G,G]$ be the \emph{commutator subgroup} or
\emph{derived group} of $G$, i.e. the subgroup generated by
commutators. $G'$ is invariant under any homomorphism, in
particular it is normal. It is the smallest normal subgroup of $G$
with an abelian factor group. Denoting $G^{(0)}:=G$, $G^{(1)}:=G'$,
$G^{(n)}:=(G^{(n-1)})'$, $n\ge 2$, one obtains \emph{derived series}
of $G$:
\begin{equation}\label{eq:derivedseries}
    G=G^{(0)}\supset G'\supset G^{(2)}\supset \dots\supset G^{(n)}\supset
\dots
\end{equation}
If $G^{(n)}=e$ for some $n$, i.e. the series (\ref{eq:derivedseries})
    stabilizes by trivial group,
    the group $G$ is \emph{solvable};

\begin{dfn}\label{dfn:polycgroup}
A solvable group with derived series with cyclic factors is called
\emph{polycyclic group}.
\end{dfn}

\begin{teo}\label{teo:almostpolaresep}
Any polycyclic-by-finite group is a strongly twisted conjugacy separable group.
\end{teo}

\begin{proof}
The class of polycyclic-by-finite groups is closed under
taking semidirect products by $\Z$. Indeed, let $G$ be an
 polycyclic-by-finite group. Then there exists a characteristic
 (polycyclic) subgroup $P$ of finite index in $G$. Hence,
 $P\rtimes \Z$ is a polycyclic normal group of $G\rtimes \Z$
 of the same finite index.

Polycyclic-by-finite groups are
conjugacy separable (\cite{Remes,Form}, see also \cite[Ch.~4]{DSegalPoly}).
It remains to apply Theorem \ref{teo:conjsepandRP}.
\end{proof}

\section{Twisted conjugacy separability and extensions}
It is known that conjugacy separability does not respect extensions.
In particular, in \cite{Goryaga} an example of a
group $G$ which is not conjugacy separable, but contains a subgroup $H$ of index
2 which is conjugacy separable, is given.

For twisted conjugacy separable groups the situation is much better under
some finiteness conditions. More precisely one has the following statement.

\begin{teo}\label{teo:phiconjandexten}
Suppose, there exists a commutative diagram
\begin{equation}\label{eq:extens}
 \xymatrix{
0\ar[r]&
H \ar[r]^i \ar[d]_{\phi'}&  G\ar[r]^p \ar[d]^{\phi} & G/H \ar[d]^{\ov{\phi}}
\ar[r]&0\\
0\ar[r]&H\ar[r]^i & G\ar[r]^p &G/H\ar[r]& 0,}
\end{equation}
where $H$ is a finitely generated normal subgroup of a finitely generated group $G$.
Suppose, $R(\phi)<\infty$, $G/H$ is a {\rm FC} group, i.e., all conjugacy
classes are finite, and $H$ is a $\phi'$-conjugacy separable group. Then
$G$ is a $\phi$-conjugacy separable group.
\end{teo}

Another variant of finiteness is $|G/H|<\infty$ (without the property $R(\phi)<\infty$).
In the relation to the next theorem, note that examples of
a conjugacy separable subgroup of
index 2 in a conjugacy non-separable group are known \cite{Goryaga}. The proof of the
next theorem is related to \cite[Prop. 3.6]{Martino}.

\begin{teo}
Let $H$ be a characteristic strongly twisted conjugacy
separable subgroup of finite index in $G$.
Then $G$ is a strongly twisted conjugacy separable group.
\end{teo}

\begin{proof}
Obviously $G$ is strongly twisted conjugacy
separable if and only if $\{g\}_\phi$ is closed in the profinite topology,
where $\phi:G\to G$ is an arbitrary automorphism and $g\in G$ is an arbitrary element.

Suppose, $x_i$, $i=1,\dots,r$, are coset representatives, where $r$ is the index of $H$
in $G$, and $g_i:=x_i g \phi(x_i^{-1})$. Then
\begin{eqnarray*}
\{g\}_\phi&=&\{w g \phi(w^{-1})\: | \: w\in G \}
=\bigcup_i \: \{hx_i g \phi(x_i^{-1})\phi(h^{-1})\: | \: h\in H \}\\
&=&\bigcup_i \: \{h g_i \phi(h^{-1})\: | \: h\in H \}=
\bigcup_i \: \{(h g_i \phi(h^{-1})g_i^{-1})g_i\: | \: h\in H \}\\
&=& \bigcup_i \: \{ e \}_{\t_{g_i}\circ \phi'}\cdot g_i.
\end{eqnarray*}
Cosets by a normal subgroup of finite index are evidently closed and open
in the profinite topology. Thus, by the supposition $\{ e \}_{\t_{g_i}\circ \phi'}$
is closed in $H$ and $G$, as well as their right translations, i.e., entries of
the above union. Since the union is finite, $\{g\}_\phi$ is closed in $G$.
\end{proof}

\section{Examples and counterexamples}\label{sec:osin}
Some of examples of groups, for which the twisted  Burnside-Frobenius
theorem in the above formulation is true, out of the class of
polycyclic-by-finite groups were obtained by F.~Indukaev \cite{IndDipl}. Namely,
it is proved that wreath products $A\wr\Z$ are RP groups,
where $A$ is a finitely
generated abelian group (these groups are residually finite).

Now let us present some counterexamples to the twisted Burnside-Frobenius
theorem in the above formulation for some discrete groups with extreme properties.
Suppose, an infinite discrete group
$G$ has a finite number of conjugacy classes.
Such examples can be found in \cite{serrtrees} (HNN-group),
\cite[p.~471]{olsh} (Ivanov group), and \cite{Osin} (Osin group).
Then evidently, the characteristic function of the unity element is not
almost-periodic and the argument above is not valid. Moreover, let us
show, that these groups give rise counterexamples to the above theorem.

In particular, they are not twisted conjugacy separable. Evidently,
they are not conjugacy separable, because they are not residually finite.

\begin{ex}\label{ex:osingroup}
For the Osin group the Reidemeister number $R(\Id)=2$,
while there is only trivial (1-dimensional) finite-dimensional
representation.
Indeed, Osin group is
an infinite finitely generated group $G$ with exactly two conjugacy classes.
All nontrivial elements of this group $G$ are conjugate. So, the group $G$
is simple, i.e. $G$ has no nontrivial normal subgroup.
This implies that group $G$ is not residually finite
(by definition of residually finite group). Hence,
it is not linear (by Mal'cev theorem \cite{malcev}, \cite[15.1.6]{Robinson})
and has no finite-dimensional irreducible unitary
representations with trivial kernel. Hence, by simplicity of $G$, it has no
finite-dimensional
irreducible unitary representation with nontrivial kernel, except of the
trivial one.

Let us remark that Osin group is non-amenable, contains the free
group in two generators $F_2$,
and has exponential growth.
\end{ex}

\begin{ex}\label{ex;ivanovgroup}
For large enough prime numbers $p$,
the first examples of finitely generated infinite periodic groups
with exactly $p$ conjugacy classes were constructed
by Ivanov as limits of hyperbolic groups (although hyperbolicity was not
used explicitly) (see \cite[Theorem 41.2]{olsh}).
Ivanov group $G$ is infinite periodic
2-generator  group, in contrast to the Osin group, which is torsion free.
The Ivanov group $G$ is also a simple group.
The proof (kindly explained to us by M. Sapir) is the following.
Denote by $a$ and $b$ the generators of $G$ described in
\cite[Theorem 41.2]{olsh}.
In the proof of Theorem 41.2 on  \cite{olsh}
it was shown that each of elements of $G$
is conjugate in $G$ to a power of generator $a$ of order $s$.
Let us consider any normal subgroup $N$ of $G$.
Suppose
$\gamma \in N$. Then $\gamma=g a^sg^{-1}$ for some $g\in G$ and some $s$.
Hence,
$a^s=g^{-1} \gamma g \in N$ and from  periodicity of $a$, it follows that also
$ a\in N$
as well as $ a^k \in N$  for any $k$, because $p$ is prime.
Then any element $h$ of $G$ also belongs to $N$
being of the form $h=\til h a^k (\til h)^{-1}$, for  some $k$, i.e., $N=G$.
Thus, the group $G$ is simple. The discussion can be completed
in the same way as in the case of Osin group.
\end{ex}

\begin{ex}
In paper \cite{HNN}, Theorem III and its corollary,
G.~Higman, B.~H.~Neumann, and H.~Neumann
proved that any locally infinite countable group $G$
can be embedded into a countable group $G^*$ in which all
elements except the unit element are conjugate to each other
(see also \cite{serrtrees}).
The discussion above related Osin group remains valid for $G^*$
groups.
\end{ex}

\section{Twisted Dehn conjugacy problem }\label{sec:desisionprob}

The subject is closely related to some decision problem. Recall that
M.~Dehn in 1912 \cite{Dehn} (see \cite[Ch. 1, \S 2; Ch. 2, \S 1]{LyndSchupp}))
has formulated in particular

\smallskip\noindent
\textbf{Conjugacy problem:} Does there exists an algorithm to determine
whether an arbitrary pair of group words $U$, $V$ in the generators of $G$ define
conjugate elements of $G$?

\smallskip
The following question was posed by G.~Makanin \cite[Question 10.26(a)]{Kourovka}:

\smallskip\noindent
\textbf{Question:} Does there exists an algorithm to determine
whether for an arbitrary pair of group words $U$ and $V$
of a free group $G$ and an arbitrary automorphism $\phi$ of $G$
the equation $\phi(X)U=VX$ solvable in $G$?

\smallskip
In  \cite{BardBokutVesnin} the following problem, which generalizes the two
above problems, was posed:

\smallskip\noindent
\textbf{Twisted conjugacy problem:} Does there exists an algorithm to determine
whether for an arbitrary pair of group words $U$ and $V$ in the generators of $G$
the equality $\phi(X)U=VX$ holds for some $W\in G$ and $\phi\in H$, where $H$ is
a fixed subset of $\Aut(G)$?

\smallskip
In  \cite{BardBokutVesnin} a partial affirmative question to the Makanin's question
is obtained.

We will discuss the twisted conjugacy problem for $H=\{\phi\}$.

\begin{teo}\label{teo:twistconjprobforpolbfin}
The twisted conjugacy problem has the affirmative answer for $G$ being
polycyclic-by-finite group and $H$ be equal to a unique automorphism $\phi$.
\end{teo}

\begin{proof}
It follows immediately from Theorem \ref{teo:almostpolaresep} by the same
argument as in the paper of
Mal'cev \cite{MalcevIvanovo} (see also \cite{Mostowski},
where the property
of conjugacy separability was first formulated)
for the (non-twisted) conjugacy problem.
\end{proof}
In fact we have proved the following statement.
\begin{teo}
If $G$ is strongly twisted conjugacy separable then the twisted
Dehn conjugacy problem is solvable for any  automorphism of $G$.
\end{teo}

Also one can study some more particular cases of this problem.
In particular, one has
\begin{teo}
Let $G$ be a $\phi$-conjugacy separable group. Then the twisted
Dehn conjugacy problem is solvable for $\phi$.
\end{teo}

The results of Section \ref{sec:resfin} imply
\begin{teo}
If $G$ is a finitely generated
residually finite group and $R(\phi)<\infty$,
then the twisted
Dehn conjugacy problem is solvable for $\phi$.
\end{teo}

From Corollary 3.4 and Proposition 3.5 in \cite{Martino} one can
obtain
\begin{teo}
Suppose $G$ is the fundamental group of a closed hyperbolic surface
and $\phi:G\to G$ is virtually inner. Then the twisted
Dehn conjugacy problem is solvable for $\phi$.
\end{teo}

\section{Some questions}\label{sec:quest}
It is evident, that any conjugacy separable group is residually
finite (because the unity element is an entire conjugacy class).
This argument does not work for general Reidemeister classes.
In this relation we wish to formulate several questions:

\smallskip
\textbf{Question 1:} Does the $\phi$-conjugacy separability imply
residually finiteness?

\smallskip
\textbf{Question 2:} Does the $\phi$-conjugacy separability imply
residually finiteness, provided $R(\phi)<\infty$?

\smallskip
\textbf{Question 3:} Does the twisted conjugacy separability imply
residually finiteness, provided the existence of $\phi$ with
$R(\phi)<\infty$?

\smallskip
\textbf{Question 4:} Let $G$ be a residually finite group and
$\phi$ its automorphism with $R(\phi)<\infty$. Is $G$
$\phi$-conjugacy separable ?

\smallskip
The affirmative answer to the last question implies twisted
Burnside-Frobenius theorem for $\phi$. This will be made in
Section \ref{sec:resfin}.

\section{Residually finite groups}\label{sec:resfin}

Recall that $G$ is called \emph{residually finite}, if for any $g\in G$
there exists an (epi)morphism $F_g$ of $G$ onto a finite group $K_g$ such that
$F_g(g)\ne e\in K_g$. In other words, $\{e\}\in G$ is closed in the profinite
topology.

The following theorem is proved in \cite{FelTroResFin}.

\begin{teo}
Let $G$ be a residually finite finitely generated group and
$\phi$ its automorphism with $R(\phi)<\infty$. Then $G$ is
$\phi$-conjugacy separable.
\end{teo}

\begin{proof}[Scetch of a proof]
It is clear, that it is sufficient to prove that
 $R(\phi)\le S_f(\phi)$.
$R(\phi)$ equals the dimension of the space of twisted invariant
elements of $\ell^\infty(G)$, i.e. functionals on $\ell^1(G)$ such
that their kernels contain the closure $K_1$ in $\ell^1(G)$ of the
space of elements of the form $b-g[b]$, $g[b](x):=b(g x
\phi(g^{-1})$.

Since $R(\phi)<\infty$, $\mathrm{codim}\, K_1 = R(\phi)$, and
$K_1$ has a Banach space complement of dimension $R(\phi)$. We can
take it in a way such that it has a base $a_i \in \C[G]$,
$i=1,\dots,R(\phi)$, i.e., all $a_i$'s have a finite support. Let
$p:G \to F=G/H$ be an epimorphism on a finite group $F$ such that
it distinguishes all elements from the union of (finite) supports
of $a_i$ and $H$ is characteristic. The image of $\ell^1(G)$ under the induced
homomorphism $p_1$ is $\ell^1(F)=\C[F]$. Also $K_1$ maps
epimorphically onto the space $K_p$ of elements
$\b-p(g)[\b]=p_1(b)-p(g)[p_1(b)]=p_1(b-g[b])$  in $\C[F]$. Thus,
$\{p_1(a_i)\}$ form a basis of a complement to $K_p$ in $\C[F]$.
Decompose this (finite dimensional)
algebra $\C[F]$ into a direct sum
of matrix algebras, i.e., decompose   the left regular
representation of $F$ into irreducible ones:
$
\l_F\cong\oplus_{i=1}^N V_i\otimes V_i^*.
$
Let $K_i$ be formed by $x-\r_i(g)[x]$ in $\End V_i$.
Since $J$ is an algebra isomorphism,
$
R(\phi)=\mathrm{codim}\, K_1= \sum_i \mathrm{codim}\, K_i.
$
The last one is 1 if $\wh \phi (\r_i)= \r_i$
and 0 otherwise.
Thus, $R(\phi)\le$ the number of finite dimensional
fixed points of $\wh\phi$.
\end{proof}

%\bibliographystyle{amsplain}
%\bibliography{felsh}

\begin{thebibliography}{10}

\bibitem{Arthur}
J.~Arthur and L.~Clozel, \emph{Simple algebras, base change, and the advanced
  theory of the trace formula}, Princeton University Press, Princeton, NJ,
  1989. \MR{90m:22041}

\bibitem{BardBokutVesnin}
Valerij Bardakov, Leonid Bokut, and Andrei Vesnin, \emph{Twisted conjugacy in
  free groups and {M}akanin's question}, Southeast Asian Bull. Math.
  \textbf{29} (2005), no.~2, 209--226, (E-print arxiv:math.GR/0401349).
  \MR{MR2217530}

\bibitem{Dehn}
M.~Dehn, \emph{\"{U}ber unendliche diskontinuierliche {G}ruppen}, Math. Ann.
  \textbf{71} (1911), no.~1, 116--144. \MR{MR1511645}

\bibitem{FelshB}
A.~Fel'shtyn, \emph{Dynamical zeta functions, {N}ielsen theory and
  {R}eidemeister torsion}, Mem. Amer. Math. Soc. \textbf{147} (2000), no.~699,
  xii+146. \MR{2001a:37031}

\bibitem{FelHill}
A.~Fel'shtyn and R.~Hill, \emph{The {R}eidemeister zeta function with
  applications to {N}ielsen theory and a connection with {R}eidemeister
  torsion}, $K$-Theory \textbf{8} (1994), no.~4, 367--393. \MR{95h:57025}

\bibitem{FelBanach}
\bysame, \emph{Dynamical zeta functions, congruences in {N}ielsen theory and
  {R}eidemeister torsion}, Nielsen theory and Reidemeister torsion (Warsaw,
  1996), Polish Acad. Sci., Warsaw, 1999, pp.~77--116. \MR{2001h:37047}

\bibitem{FelHillWong}
A.~Fel'shtyn, R.~Hill, and P.~Wong, \emph{Reidemeister numbers of equivariant
  maps}, Topology Appl. \textbf{67} (1995), no.~2, 119--131. \MR{MR1362078
  (96j:58139)}

\bibitem{FelIndTro}
A.~Fel'shtyn, F.~Indukaev, and E.~Troitsky, \emph{Twisted {B}urnside theorem
  for two-step torsion-free nilpotent groups}, {C}*-algebras and elliptic
  theory. {II}, Trends in Math., Birkh\"auser, 2008, pp.~87--101.

\bibitem{FelTro}
A.~Fel'shtyn and E.~Troitsky, \emph{A twisted {B}urnside theorem for countable
  groups and {R}eidemeister numbers}, Noncommutative Geometry and Number Theory
  (C.~Consani and M.~Marcolli, eds.), Vieweg, Braunschweig, 2006, pp.~141--154.

\bibitem{polyc}
A.~Fel{\cprime}shtyn and E.~Troitsky, \emph{Twisted {B}urnside-{F}robenius
  theory for discrete groups}, J. Reine Angew. Math. \textbf{613} (2007),
  193--210. \MR{MR2377135}

\bibitem{FelTroVer}
A.~Fel'shtyn, E.~Troitsky, and A.~Vershik, \emph{Twisted {B}urnside theorem for
  type {II}${}_1$ groups: an example}, Math. Res. Lett. \textbf{13} (2006),
  no.~5, 719--728.

\bibitem{FelTroResFin}
Alexander Fel'shtyn and Evgenij Troitsky, \emph{Twisted conjugacy classes in
  residually finite groups}, Arxiv e-print 1204.3175, 2012.

\bibitem{Form}
E.~Formanek, \emph{Conjugacy separability in polycyclic groups}, J. Algebra
  \textbf{42} (1976), 1--10.

\bibitem{gowon}
D.~Gon{\c{c}}alves and P.~Wong, \emph{Twisted conjugacy classes in exponential
  growth groups}, Bull. London Math. Soc. \textbf{35} (2003), no.~2, 261--268.
  \MR{2003j:20054}

\bibitem{Goryaga}
A.~V. Goryaga, \emph{Example of a finite extension of an {FAC}-group that is
  not an {FAC}-group}, Sibirsk. Mat. Zh. \textbf{27} (1986), no.~3, 203--205,
  225. \MR{MR853899 (87m:20090)}

\bibitem{GrFA}
R.I. Grigorchuk, \emph{On {B}urnside's problem on periodic groups}, Funct.
  Anal. Appl. \textbf{14} (1980), 41--43.

\bibitem{Groth}
A.~Grothendieck, \emph{Formules de {N}ielsen-{W}ecken et de {L}efschetz en
  g\'eom\'etrie alg\'ebrique}, S\'eminaire de G\'eom\'etrie Alg\'ebrique du
  {B}ois-{M}arie 1965-66. {SGA} 5, Lecture Notes in Math., vol. 569,
  Springer-Verlag, Berlin, 1977, pp.~407--441.

\bibitem{GuSi}
N.~Gupta and S.~Sidki, \emph{On the {B}urnside problem for periodic groups},
  Math. Z. \textbf{182} (1983), 385--388.

\bibitem{HNN}
Graham Higman, B.~H. Neumann, and Hanna Neumann, \emph{Embedding theorems for
  groups}, J. London Math. Soc. \textbf{24} (1949), 247--254. \MR{MR0032641
  (11,322d)}

\bibitem{IndDipl}
Fedor Indukaev, \emph{Twisted burnside theory for infinite groups: some
  examples}, Master thesis, Moscow State University, May 2006.

\bibitem{Jiang}
B.~Jiang, \emph{Lectures on {N}ielsen fixed point theory}, Contemp. Math.,
  vol.~14, Amer. Math. Soc., Providence, RI, 1983.

\bibitem{Jiang88}
Bo~ju Jiang, \emph{A characterization of fixed point classes}, Fixed point
  theory and its applications (Berkeley, CA, 1986), Contemporary Math., no.~72,
  Amer. Math. Soc., 1988, pp.~157--160.

\bibitem{Kirillov}
A.~A. Kirillov, \emph{Elements of the theory of representations},
  Springer-Verlag, Berlin Heidelberg New York, 1976.

\bibitem{LyndSchupp}
Roger~C. Lyndon and Paul~E. Schupp, \emph{Combinatorial group theory},
  Springer-Verlag, Berlin, 1977, Ergebnisse der Mathematik und ihrer
  Grenzgebiete, Band 89. \MR{MR0577064 (58 \#28182)}

\bibitem{malcev}
A.~I. Mal'cev, \emph{On the faithful representations of infinite groups by
  matrices}, Mat. Sb. (NS) \textbf{8(50)} (1940), 405--422, (in Russian.
  English translation: Amer. Math. Soc. Transl. (2), \textbf{45} (1965),
  1--18).

\bibitem{MalcevIvanovo}
A.I. Malcev, \emph{On homomorphisms onto finite groups}, Uchen. Zapiski
  Ivanovsk. ped. instituta \textbf{18} (1958), no.~5, 49--60, (=Selected
  Papers, Vol.~1, 1976, 450--461).

\bibitem{Martino}
Armando Martino, \emph{A proof that all {S}eifert $3$-manifold groups and all
  virtual surface groups are conjugacy separable}, E-print
  arxiv:math.GR/0505565, 2005.

\bibitem{Kourovka}
V.~D. Mazurov and E.~I. Khukhro (eds.), \emph{The {K}ourovka notebook},
  augmented ed., Rossi\u\i skaya Akademiya Nauk Sibirskoe Otdelenie, Institut
  Matematiki im. S. L. Soboleva, Novosibirsk, 2002, Unsolved problems in group
  theory. \MR{MR1956290 (2004f:20001a)}

\bibitem{Mostowski}
A.~W. Mostowski, \emph{On the decidability of some problems in special classes
  of groups}, Fund. Math. \textbf{59} (1966), 123--135. \MR{MR0224693 (37
  \#292)}

\bibitem{olsh}
A.~Yu. Ol{\cprime}shanski{\u\i}, \emph{Geometry of defining relations in
  groups}, Mathematics and its Applications (Soviet Series), vol.~70, Kluwer
  Academic Publishers Group, Dordrecht, 1991, Translated from the 1989 Russian
  original by Yu.\ A. Bakhturin. \MR{MR1191619 (93g:20071)}

\bibitem{Osin}
Denis Osin, \emph{Small cancellations over relatively hyperbolic groups and
  embedding theorems}, Ann. of Math. (2) \textbf{172} (2010), no.~1, 1--39.
  \MR{2680416 (2012a:20068)}

\bibitem{Remes}
V.~N. Remeslennikov, \emph{Conjugacy in polycyclic groups. ({R}ussian)},
  Algebra i Logika \textbf{8} (1969), 712--725.

\bibitem{Robinson}
Derek J.~S. Robinson, \emph{A course in the theory of groups}, second ed.,
  Graduate Texts in Mathematics, vol.~80, Springer-Verlag, New York, 1996.
  \MR{MR1357169 (96f:20001)}

\bibitem{DSegalPoly}
Daniel Segal, \emph{Polycyclic groups}, Cambridge Tracts in Mathematics,
  no.~82, Cambridge University Press, Cambridge, 1983.

\bibitem{serrerepr}
Jean-Pierre Serre, \emph{Linear representations of finite groups},
  Springer-Verlag, New York, 1977, Translated from the second French edition by
  Leonard L. Scott, Graduate Texts in Mathematics, Vol. 42. \MR{MR0450380 (56
  \#8675)}

\bibitem{serrtrees}
\bysame, \emph{Trees}, Springer Monographs in Mathematics, Springer-Verlag,
  Berlin, 2003, Translated from the French original by John Stillwell,
  Corrected 2nd printing of the 1980 English translation. \MR{MR1954121
  (2003m:20032)}

\bibitem{Shokra}
Salahoddin Shokranian, \emph{The {S}elberg-{A}rthur trace formula}, Lecture
  Notes in Mathematics, vol. 1503, Springer-Verlag, Berlin, 1992, Based on
  lectures by James Arthur. \MR{MR1176101 (93j:11029)}

\bibitem{WiZa}
J.~S. Wilson and P.~A. Zalesskii, \emph{Conjugacy separability of certain
  torsion groups}, Arch. Math. (Basel) \textbf{68} (1997), no.~6, 441--449.

\end{thebibliography}
%\end{document}

\def\cprime{$'$} \def\cprime{$'$} \def\cprime{$'$} \def\cprime{$'$}
  \def\cprime{$'$} \def\cprime{$'$} \def\cprime{$'$} \def\dbar{\leavevmode\hbox
  to 0pt{\hskip.2ex \accent"16\hss}d} \def\cprime{$'$} \def\cprime{$'$}
  \def\polhk#1{\setbox0=\hbox{#1}{\ooalign{\hidewidth
  \lower1.5ex\hbox{`}\hidewidth\crcr\unhbox0}}} \def\cprime{$'$}
  \def\cprime{$'$} \def\cprime{$'$} \def\cprime{$'$}
\providecommand{\bysame}{\leavevmode\hbox to3em{\hrulefill}\thinspace}
\providecommand{\MR}{\relax\ifhmode\unskip\space\fi MR }
\providecommand{\MRhref}[2]{%
  \href{http://www.ams.org/mathscinet-getitem?mr=#1}{#2}
}
\providecommand{\href}[2]{#2}

\end{document}